\documentclass{amsart}[11]

\usepackage{amssymb}
\usepackage{amsmath}
\usepackage{amsthm}
\usepackage{amsbsy}
\usepackage[dvips]{graphics}

\theoremstyle{plain}
\newtheorem{theorem}{Theorem}[section]

\newtheorem{lemma}[theorem]{Lemma}

\newtheorem{proposition}[theorem]{Proposition}

\newtheorem*{algorithm}{Algorithm}

\theoremstyle{definition}
\newtheorem{definition}{Definition}[section]

\theoremstyle{remark}
\newtheorem{remark}[theorem]{Remark}
\newtheorem*{acknowledgements}{Acknowledgements}
\newtheorem{example}{Example}[section]

\def \del{\partial}

\begin{document}

\title{On the Andrews-Curtis conjecture and Algorithms from Topology}

\begin{abstract}
We relate the Andrews-Curtis conjecture to the triviality problem for
balanced presentations using algorithms from $3$-manifold
topology. Implementing this algorithm could lead to counterexamples to
the Andrews-Curtis conjecture.
\end{abstract}


\date{First corrected version : January 2000 ; this revision: August 7, 2001}

\author{Siddhartha Gadgil}

\address{	Department of mathematics\\
		253-37, California Institute of Technology\\
		Pasadena, CA 91125}
\email{gadgil@cco.caltech.edu}

\maketitle 

The Andrews-Curtis conjecture \cite {AC} says that a balanced
presentation of the trivial group can be transformed to a standard
presentation
$\left<\alpha_1,\dots,\alpha_n;\alpha_1,\dots,\alpha_n\right>$ by
using `Andrews-Curtis moves', i.e. multiplying one relation by
another, inverting a relation and conjugating a relation by a
generator. Here, we use algorithmic methods from $3$-manifold topology
to relate this conjecture to a fundamental algorithmic question in
group theory, namely the triviality problem for balanced
presentations.

Thus, our main result is

\begin{theorem}[\ref{T:bal}] At least one of the following holds
\begin{itemize}
\item There is an algorithm to recognise balanced presentations of the
trivial group, or
\item The (balanced) Andrews-Curtis conjecture is false.
\end{itemize}
\end{theorem}

Moreover, our methods could, after extensive computations, show the
failure of the (balanced) Andrews-Curtis conjecture for balanced
presentations, and also `Property-R' for links.

Our approach is to construct an algorithm, using algorithms from
$3$-manifold topology, in particular the Rubinstein-Thompson algorithm
(\cite{Ru},\cite{Th} and
\cite{Ma}) for recognising the $3$-sphere, to recognise certain
presentations of the trivial group. To do this, we associate with a
group presentation a `handle-structure', a natural generalisation of a
handle-decomposition of a $3$-manifold. Just as the
Rubinstein-Thompson-Matveev recognises a certain class of
handle-decompositions of homotopy $3$-spheres, namely those
corresponding to a $3$-sphere, our algorithm recognises
handle-structures corresponding to certain presentations of the
trivial group. It is not clear exactly what this class of
presentations is, but one can see that it includes those related to
the standard presentations by the Andrews-Curtis moves.

We can immediately conclude that either the Andrews-Curtis conjecture
is false, or we have an algorithm to solve the triviality problem for
balanced presentations. Indeed, if the Andrews-Curtis conjecture were
true, then our algorithm would recognise the trivial group given a
balanced presentation.

Further, we can apply the algorithm to candidate counterexamples. If
the answer is in the negative, then we know that we have a
presentation that violates the Andrews-Curtis conjecture. In
particular, if such a presentation comes from a Kirby diagram for
$S^3$ then property-R for links is false.

\begin{acknowledgements}
	I would like to thank David Gabai for his indispensable help
and encouragement at all stages of this work. I am indebted to Michael
Freedman, Andrew Casson, Robion Kirby, Peter Teichner, Dave Letscher
and Tao Li for many valuable comments and conversations.
\end{acknowledgements}

\section{Handle-structures}\label{S:hndl}

A connected $3$-manifold with non-empty boundary has a
handle-decomposition with a $0$-handle, some $1$-handles and some
$2$-handles. We call $0$, $1$ and $2$-handles balls beams and plates
respectively. A beam is homeomorphic to $D^2\times [0,1]$ and is
attached to the ball along $D^2\times \{0,1\}$. We shall call
$D^2\times \{0,1\}$ the `sticky end' of a beam. Likewise a plate is
homeomorphic to $D^2\times [0,1]$ and is attached to the union of
balls and beams along $S^1\times [0,1]$, which we shall call the
`sticky end' of the plate. The components of intersection of beams and
plates with the ball shall be called islands and bridges respectively.
 
It is clear that a handle-decomposition as above is completely
determined by the beams and plates together with their `sticky ends',
for there is an essentially unique way to stick the `sticky ends' to
the $0$-handle. Further, suppose we are given beams and plates, with
the plates glued to beams along strips in their sticky end as they
would be in a $3$-manifold. Then these are the beams and plates of a
handle-decomposition of a $3$-manifold iff the sticky end of these is
planar. This is an unnecessary assumption for many applications of
$3$-manifold techniques, in particular many results in normal surface
theory. This motivates the definition

\begin{figure}
\includegraphics{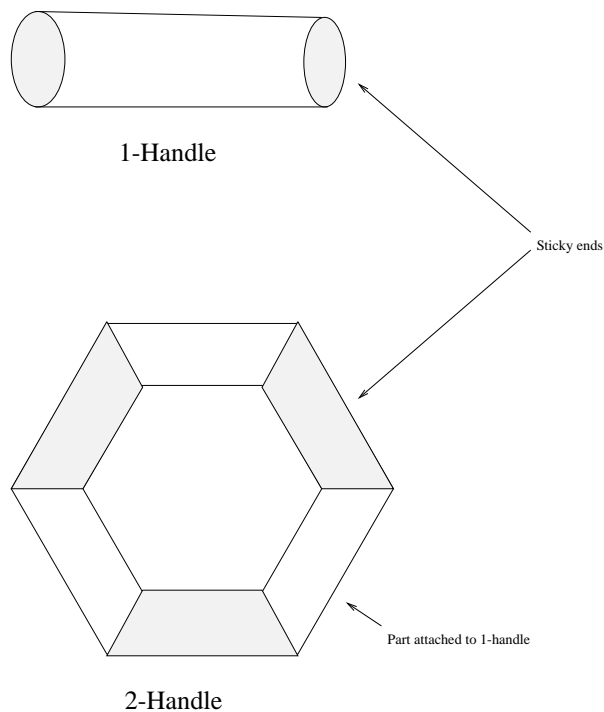}
\caption{$1$-handles and $2$-handles in a handle structure}\label{F:sticky}
\end{figure}

\begin{definition} A {\em handle-structure} is a union of beams and
plates, which intersect in unions of discs that when viewed as subsets
in either the beams or the plates are of the form $\alpha\times [0,1]$
where $\alpha$ is an arc in $D^2$. The horizontal boundaries
$D^2\times \{0,1\}$ of the beams together with the portion of the
horizontal boundary of the plates that is not glued to a beam form a
surface consisting of islands and bridges that we call the {\em sticky
end} (see figure~\ref{F:sticky}). The space consisting of the beams
and plates will be called the {\em total space}.
\end{definition}

Our aim in introducing handle-structures is to associate these with
finite presentations of groups. We shall use the notation $X/Y$, where
$X$ and $Y$ are topological spaces with $Y\subset X$, to denote the
quotient space obtained form $X$ by identifying $Y$ to a
point. Briefly, we can associate with any presentation a singular
$3$-manifold, with the only singularities being transversal
intersections and self-intersections of the bridges. We delete from
this the closure of the ball to get a non-singular open
$3$-manifold. The completion of this gives a handle-structure, with
the sticky-part being the points added to complete the space.

\begin{proposition}\label{P:gphan} Given a finite presentation of a
group $G$, we can associate to it a handle-structure. This
construction is canonical up to finitely many choices. Further, if $M$
is the total space and $S$ the sticky end, then $\pi_1(M/S)=G$.
\end{proposition}
\begin{proof} Take a beam corresponding to each generator and a plate
corresponding to each relation. Successively glue the plates to the
beams according to the word that the relation represents. Namely,
start with a base point that will not be glued to any $1$-handle, and
for each successive letter in the word of the relation attach a strip
in the horizontal boundary $S^1\times [0,1]$ of the plate to the
corresponding beam, with the direction determined by the sign of the
exponent. Thus, strips attached to $1$-handles alternate with bridges
as in figure~\ref{F:sticky}, with the attaching maps being governed by
the word that the $2$-handle represents. The only essential choices
involved are in choosing where to attach the strip to the beam
relative to previous intersections with plates.

Clearly $M/S$ deformation retracts onto a $2$-complex that is well
known to have fundamental group $G$. More precisely, each one handle
deformation retracts to its core, and in the quotient the two islands
of a $1$-handle are identified to a single point. Thus, the
$1$-handles deformation retract onto a wedge of circles. The
$2$-handles in turn deformation retract on to their cores, which are
$2$-discs. By construction, the attaching maps of the disc
corresponding to a given $2$-handle to the wedge of circles formed by
the $1$-handles gives the corresponding relation in the fundamental
group.
\end{proof}

\begin{figure}
\includegraphics{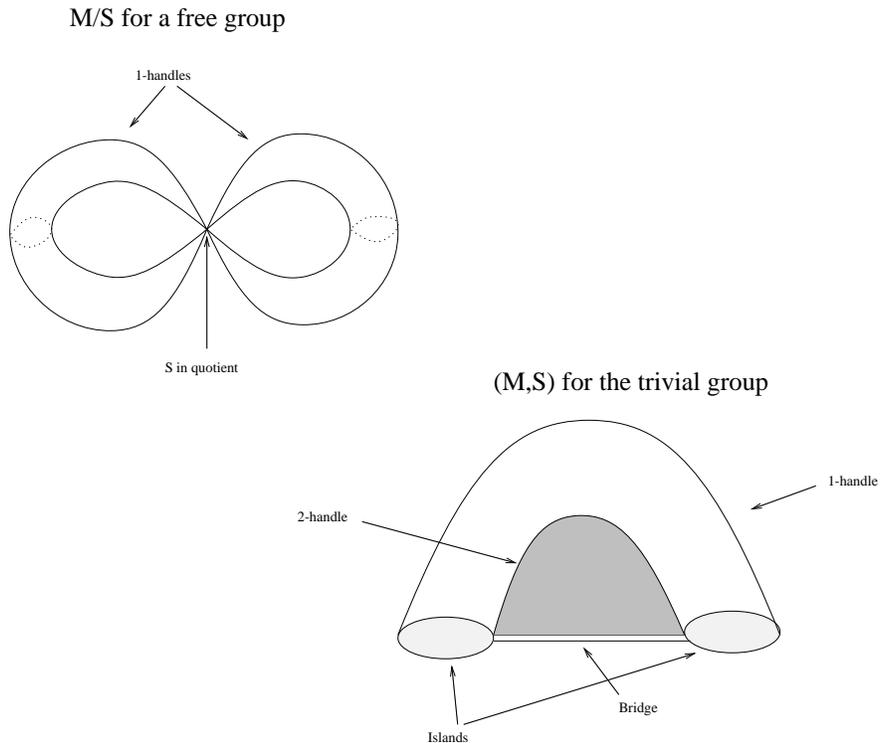}
\caption{Examples of $M$ and $S$}\label{F:eg}
\end{figure}

\begin{example} 
For the presentation $\left<\alpha_1,...,\alpha_n;\right>$ of a
free group on $n$ generators, we have $n$ $1$-handles with $2n$
islands, $2$ on each $1$-handle, and no $2$-handles. Thus $M/S$ is a
thickened wedge of $n$ circles (see figure~\ref{F:eg}).
\end{example}

\begin{example} 
For the presentation $\left<\alpha;\alpha\right>$ of the trivial
group, we have one $1$-handle and a $2$-handle, with the intersection
of these connected. Thus $M$ is the union of two $3$-balls glues
together along a $2$-disc, and is therefore a $3$-ball. $S$ consists
of two islands joined together by a bridge, and hence is $B^3$ (see
figure~\ref{F:eg}).
\end{example}

\begin{example} The standard balanced presentation of the trivial group
with $n$ generators and relations gives a handle-structure with $n$
components, each of them as in the previous example.
\end{example}

We shall recognise handle-structures that correspond to a class of
presentations of the trivial group.

\section{The Algorithm}\label{S:alg}

We assume for the next two sections that we are given a {\em balanced}
presentation of a group. To this we associate a finite collection
of handle structures as above. We apply the following algorithm to
each of these.

\begin{algorithm} Given a handle-structure, let $M$ be the total space
and $S$ the sticky end.
\begin{enumerate}

\item\label{loop} Attach thickened discs to $M$ along each component of
$\del S$ to get a closed manifold. Define the sticky part of this to
be the union of the components of the boundary that intersect $S$.

\item\label{RTM} If each component of the resulting manifold is
$S^2\times [0,1]$, with one of its boundary components being the
sticky part, then answer `Yes' and terminate.

\item\label{Ann} Cut $(M,S)$ along a maximal family of pairwise
disjoint and non-parallel annuli such that each annulus has exactly
one boundary component in $S$. If there are no such annuli, answer `No'
and terminate.

\item Repeat from step \ref{loop} using the new $(M,S)$.\qed

\end{enumerate}
\end{algorithm}

We shall call annuli which are as in step~\ref{Ann} of the algorithm
\emph{vertical annuli}.

\begin{remark} Step \ref{RTM} can be achieved by using the
Rubinstein-Thompson-Matveev algorithm to recognise $S^3$. Step
\ref{Ann} uses basic normal surface theory. The algorithm terminates
because the boundary gets simpler at each iteration.
\end{remark}

\begin{lemma}\label{L:gptriv} If the algorithm terminates with a `Yes',
then the $\pi_1(M/S)$ is the trivial group.
\end{lemma}
\begin{proof} This follows by induction on the number of steps
needed. If $M$ is as in step~\ref{RTM} then this is true. It is easy
to see that this remains true if $S^2\times [0,1]$ is obtained after
cutting along vertical annuli and filling in discs as in
steps~\ref{Ann}~\&~\ref{loop}. We need to see that if $\pi_1(M/S)$ is
trivial after cutting along an annulus or capping off, then it is
trivial before making these changes. In the latter case, this follows
as there is no relation added - the boundary of the disc is already
trivial as it is in $S$. In the case where we cut along a vertical
annulus, let $M'$ and $S'$ be obtained after cutting. Now an arc that
in $M$ that intersects the annulus can be homotoped in the quotient
off this - it breaks into an arc ending in $S$ just before the annulus
and an arc starting from S just after the annulus. This is in $M'$, hence
trivial in $M'/S'$, and so all the more in $M/S$.
\end{proof}

We thus have an algorithm to recognise a certain class of
presentations of the trivial group. Namely, given a presentation, we
apply the above algorithm to all the handle-structures associated to
it. If the algorithm says `Yes' for at least one of these, we conclude
that the presentation corresponds to the trivial group. We shall call
presentations where this happens, and also the handle-structures for
which the algorithm says `Yes' {\em sphere-like}.

\begin{remark} In case the presentation and the handle-structure come
from a $3$-manifold, after capping off $S$ is a union of $2$-spheres,
and handle-structure is sphere-like iff the manifold is a sphere.
\end{remark}

\section{Relation to the Andrews-Curtis moves}\label{S:ac}

We have constructed above an algorithm to recognise certain
presentations of the trivial group, which we call {\em
sphere-like}. For this to be useful, we need to see that this class of
presentations is reasonably large. In fact, we shall see

\begin{theorem}\label{T:AC} A balanced presentation of the trivial
group that comes from the standard one by Andrews-Curtis moves is
sphere-like.
\end{theorem}
\begin{proof} We shall show this in two stages. First, we multiply
relations merely by concatenation without any cancellation. We shall
make a choice of handle-structures that correspond to the new (in
general not reduced) presentation such that the new presentation is
sphere-like. We then consider the effect of cancellation and show that
even after cancellation, the handle-structure is sphere-like.

The first step is accomplished by induction on the number of moves. In the
case of the standard presentations, $(M,S)=\cup (D^2\times [0,1],D^2
\times 0)$, which is sphere-like. 

We shall have to consider two basic transformation of the pairs
$(M,S)$ --  \emph{gluing} and \emph{puncturing}.

\begin{definition} Let $x$ and $y$ be two points on $\del S$. Then by
\emph{gluing} $x$ and $y$ we mean the identification of disjoint
regular neighbourhoods of these points, with the new $S$ being the
image of $S$ under the quotient map.
\end{definition}

Note that the new total space and the new sticky part are obtained
from the old ones by taking a $\del$-connected sum.

\begin{definition} Let $\gamma$ be an arc properly embedded in $M$ with
exactly one end-point in $S$. By \emph{puncturing} along $\gamma$ we
mean deleting a regular neighbourhood of $\gamma$.
\end{definition}

The arc above need not be unknotted. We shall call such an arc a
\emph{vertical arc}.

\begin{lemma} On making an appropriate choice of handle-structure
corresponding to a presentation obtained by multiplying words without
cancellation, the effect of the Andrews-Curtis moves is a sequence of
gluings and puncturings.
\end{lemma}

\begin{proof} Firstly, the inversion of a relation has no effect on
the associated handle-structures (with obvious choices). 

Next consider conjugation of a relation by a generator. This, again
with an obvious choice of handle structure, amounts to a handle slide
together with puncturing a hole along an arc with one boundary point
on the sticky part. It may be helpful to consider the change in handle
structure on making the Andrews-Curtis move
$\left<\alpha,\beta;\alpha,\beta\right>\to\left<\alpha,\beta;\alpha,
\alpha^{-1}\beta\alpha\right>$. The resulting $(M,S$ is
$(Annulus\times [0,1], Annulus)$ (see figure~\ref{F:eg3})). The
transformation in general is very similar to this particular case. If
we attach a thick disc along the new bridge that forms a loop, we see
that the resulting $(M,S)$ is obtained from the original one by
gluing. Puncturing along the arc dual to the thick disc now gives us
the handle-structure obtained by conjugation.

Next, we consider the effect of multiplying a relation by
another. First consider the particular example
$\left<\alpha,\beta;\alpha,\beta\right>
\to\left<\alpha,\beta;\alpha,\alpha\beta\right>$. The new
handle-structure is $(M,S)=(D^3,D^2)$ (see figure~\ref{F:eg2}) and is
clearly obtained by a gluing. Suppose that we are replacing the
relation $r_1$ by $r_1r_2$. If $r_2$ has only one bridge, then we have
a similar situation. Namely, in the new handle-structure, $r_1$ has a
bridge replaced by two bridges and has an additional component of
intersection with the only beam that $r_2$ intersects. Choose the new
handle structure so that this component of intersection is adjacent to
the attachment of $r_2$ to the beam. Then the handle-structure is
clearly obtained topologically by gluing.

Suppose now $r_2$ has more than one bridge, i.e., bridges in addition
to the one along which gluing has taken place.. We attach the
relations other than $r_1$ as before. The new $r_1$ is the old one
with additional bridges and gluings to $1$-handles. Namely first pick
bridges of $r_2$ and $r_1$ corresponding to the base points, and
replace the bridge of $r_1$ by two bridges to the islands bounding the
bridge of $r_2$ and adjacent to this bridge . Then attach $r_1$ to the
$1$-handles adjacent to the where $r_2$ has been attached. The new
bridges of $r_1$ will thus be adjacent to those of $r_2$. Attaching
thick discs to these gives, as in the example above, a
handle-structure obtained by gluing. As in the case of conjugation,
the final structure is obtained by puncturing along dual arcs to these
discs.
\end{proof}

\begin{figure}
\includegraphics{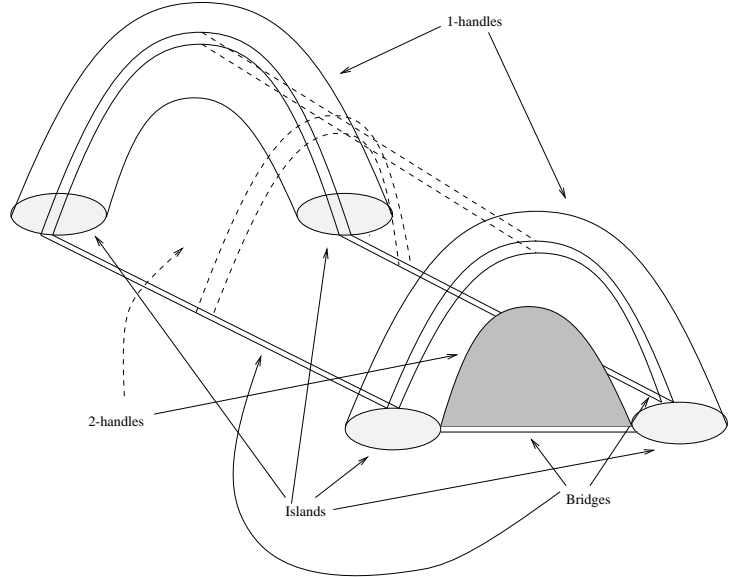}
\caption{The presentation $<\alpha,\beta;\alpha,\alpha\beta>$}\label{F:eg2}
\end{figure}

\begin{figure}
\includegraphics{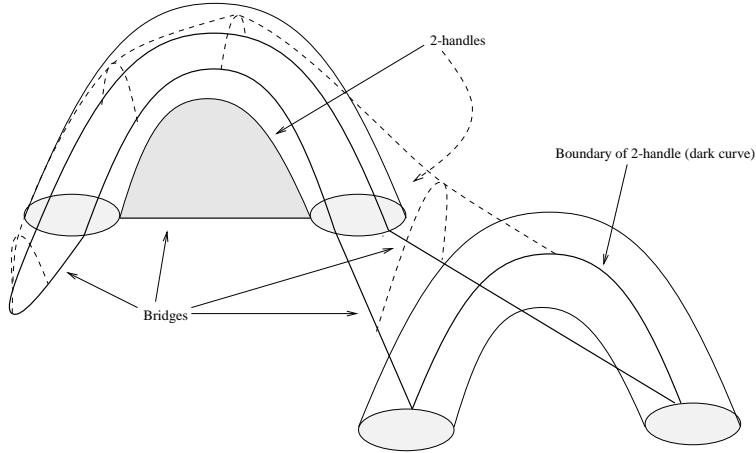}
\caption{The presentation
$<\alpha,\beta;\alpha,\alpha^{-1}\beta\alpha>$}\label{F:eg3}
\end{figure}

We shall need some further transformations to take into account
cancellations.

\begin{definition} A \emph{bigon} in $(M,S)$ is a properly embedded
disc in $M$ whose boundary has exactly one component of intersection
with each of $S$ and $\del M\setminus S$.
\end{definition}

\begin{remark} The image of the discs identified in gluing forms a
bigon in the resulting manifold. Thus, we shall sometimes refer to
gluing as gluing along a bigon.
\end{remark}

\begin{definition} A  \emph{quadrilateral} in $(M,S)$ is a properly
embedded disc in $M$ whose boundary has exactly two components of
intersection with each of $S$ and $\del M\setminus S$. 
\end{definition}

\begin{lemma} The handle-structure obtained after cancelling a pair of
adjacent letters that are inverses of each other is obtained by cutting
along a quadrilateral and then along a bigon.  
\end{lemma}

\begin{proof} To see the above, note that when we have a cancelling
pair of letters, the $2$-handle goes over a $1$-handle and then
returns over the same $1$-handle immediately. Thus, succesive bridges
of the $2$-handle are one ending in the $1$-handle, one beginning and
ending in the $1$-handle (which shall often be referred to as a
\emph{loop} in what follows) and then one beginning in the
$1$-handle. We cut along a quadrilateral (see figure~\ref{F:quad})
consisting of separating arcs along the first and last of these
bridges and arcs in the vertical boundary of the $2$-handle. Then we
cut along a bigon whose boundary has an arc in the bridge that begins
and ends in the $1$-handle. The result is homeomorphic to a
handle-structure corresponding to the cancelled presentation. Note
that in the figure, after splitting along the quadrilateral, we
declare the left side, i.e., the side away form the bigon as part of
the sticky end. Thus what we have left is indeed a bigon.

\begin{figure}
\includegraphics{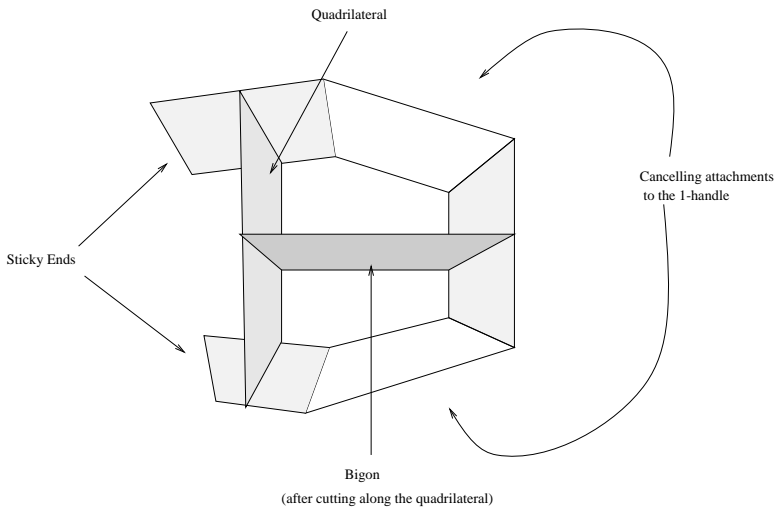}
\caption{The quadrilateral and the bigon}\label{F:quad}
\end{figure}

It may be helpful here to see these in figure~\ref{F:eg3}. Here, The
quadrilateral splits the two bridges joining the two $1$-handles. It
is thus the dotted arc joining these bridges in the figure
thickened. The side towards the lower handle in the figure becomes
part of $S$. The bigon is the dotted line joining the third bridge to
the quadrilateral thickened.
\end{proof}

We have thus reduced the problem to showing that puncturing, gluing
and cutting along bigons and quadrilaterals coming from cancellations
takes sphere-like pairs $(M,S)$ to sphere-like pairs.

\begin{lemma} Gluing and puncturing take sphere-like pairs $(M,S)$ to
sphere-like pairs.
\end{lemma}

\begin{proof}  The effect of puncturing along a vertical arc is
undone in the first step of the algorithm, where the components of $S$
are capped off by attaching discs. Thus, there is nothing to prove in
this case.

When we glue along a bigon, there are two cases to consider, namely
where we glue two distinct boundary components and where we glue a
boundary component to itself. In the latter case, the manifold after
capping off does not change, we merely glue two discs instead of
one. In the former case, after gluing and capping off, we end up with
an essential vertical annulus. This is because gluing together and
capping off has the same effect as attaching a cylinder between the
two boundary components. After cutting along this and capping off, we
end up with the manifold obtained by capping off the original
annulus. Finally, we note that as though the maximal family of
vertical annuli need not be unique, the complementary region does not
depend on the family chosen.  Hence we may assume that the family
chosen in the algorithm contains the abovementioned annulus.
\end{proof}

\begin{lemma} Cutting along a bigon takes a sphere-like pair $(M,S)$ to
a sphere-like pair.
\end{lemma}

\begin{proof} Note that there are at most two punctures intersecting
the bigon, corresponding to the two `vertices' of the bigon. If the
vertices are on different components of $\del S$, then we see that
cutting along the bigon and capping off has the same effect as capping
off the original manifold.

Thus, we need to consider only the case where both vertices are on the
same component. We first consider the case where we have a punctured
product $(\tilde S\times I)\setminus punctures$, with the only
punctures being ones that intersect the bigon. We shall then reduce to
this case. Note that we can assume that punctures disjoint from the
bigon were made after cutting along the bigon and thus the latter part
of the reduction is immediate.

Thus, the vertices are either on a puncture or a component $\gamma$ of
$\del\tilde S$. In the former case, we call a neighbourhood of the
puncture the \emph{boundary annulus} and in the latter case we use
this term for $\gamma\times I$. We shall see that we have two
possibilities, a bigon parallel to an (in general) knotted puncture
and a possibly essential one in a product (i.e. with vertex on
$\del\tilde S$ or an unknotted puncture). 

We have a natural vertical annulus in the manifold capped off at the
component of $\del S$ containing the vertex (see
figure~\ref{F:bigon}).  Namely, we add to the bigon the neighbourhood
of an arc in $\del S$, i.e., a strip in the boundary annulus. If this
is an essential vertical annulus, then as the arc along which the
manifold was punctured is a vertical arc on this annulus, puncturing
and then cutting along the bigon still gives a product. Otherwise, the
boundary curves of the annulus must bound discs, giving a
$2$-sphere. This must bound a $3$-ball (which follows from
irreducibility of the punctured product). Thus, cutting along the
bigon gives one component $D^2\times I$ and the other homeomorphic to
the original manifold. In either case, the manifold is still
sphere-like.

\begin{figure}
\includegraphics{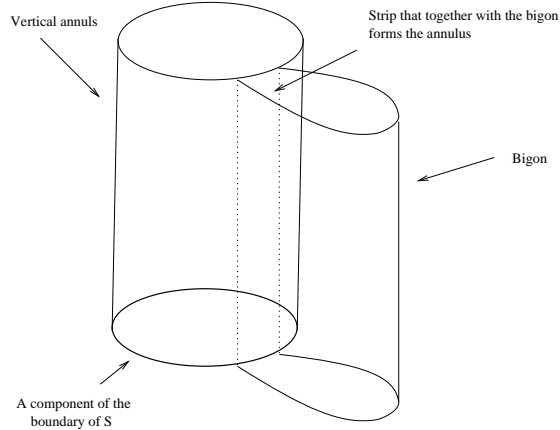}
\caption{Splitting along a bigon}\label{F:bigon}
\end{figure}

Now to reduce to the case considered, we need to consider the
intersection of the annulus with the bigons along which gluing has
taken place. This intersection, after making things in general
position, is a collection of inessential arcs (i.e. those with both
endpoints on the same boundary component of the annulus), vertical
arcs and horizontal circles. We take an outermost such component. If
it is an inessential arc, then we can isotope the annulus to remove
this intersection. To do this, note that the arc separates a disc in
each of the annulus and the bigon, and these together form a
compressing disc. As the boundary is incompressible, we have a sphere,
which must bound a ball by irreducibility. If we have a horizontal
circle, then this bounds a disc in the annulus, and we see that we
have a bigon parallel to a puncture.

Finally, in case of a vertical arc, our bigon decomposes into a piece
between two parallel punctures (perhaps the same one and a bigon with
simpler intersection. Again, the first bigon is as in one of the above
cases, and in either case, it is easy to see that we get a sphere-like
pair.

\begin{lemma}\label{T:quad} Cutting along a quadrilateral corresponding to a
cancellation takes a sphere-like pairs $(M,S)$ to a sphere-like pair.
\end{lemma}

First, observe that when cutting along a quadrilateral, there is a
choice involved. Namely, the disc along which $M$ has been cut now
corresponds to two discs in $\del M$, and one of these is to be
included in $S$. If the wrong choice is made, the result is not even
homologically $S\times I$. It is easy to see that in our case we do
have a product homologically, and so the choice has been made
correctly. Namely, if we consider the handle-structure obtained with
the other choice, then we have two bridges, one corresponding to the
earlier loop and one a new one, so that the bridges end on both sides
in the two islands corresponding to the one-handle where the
cancellation takes place. For each bridge, we can take an arc along it
and close it up in the one handle to get a closed curve. These curves
bound an annulus in $M$, and are thus homologically equivalent. But an
arc across one of the bridges intersects one of the curves once, but
not the other, hence the curves are not homologous in $S$. Thus we do
not have a (homological) product with this wrong choice.

\begin{figure}
\includegraphics{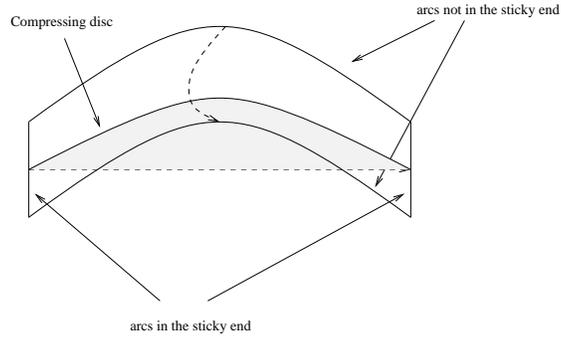}
\caption{Compressing the quadrilateral}\label{F:compress}
\end{figure}

With the correct choice, this claim is well known in the case of
products (and more generally for taut sutured manifolds~\cite{Ga}). In
our situation, we claim that cutting along a quadrilateral has the
effect of cutting along two bigons and then gluing along a bigon. More
specifically, we can find $\del$-compressing disc for the
quadrilateral, i.e., a disc whose boundary consists of an arc in the
quadrilateral and an arc in $S$ (see figure~\ref{F:compress}) at some
stage in the algorithm. Boundary-compressing along such a disc gives
two bigons, and it is easy to see that cutting along the quadrilateral
is equivalent to cutting along each of these discs and then re-gluing
along the $\del$-compressing disc.

We show that such a disc exists by induction. In the case of a
product, this follows by cutting up $M$ along vertical annuli and
discs. An outermost disc of the vertical discs and annuli cut by the
quadrilaterals gives the compression, except in the trivial case where
the quadrilateral can be made disjoint from these. Similarly, if we
have a quadrilateral after gluing along bigons, we look at outermost
arcs of the bigons after cutting by the quadrilateral to reduce to the
situation before gluing along the bigon.
        
As with bigons, we may assume that punctures not involving the
vertices of the quadrilateral were made after cutting along the
quadrilateral. Further, we may again assume that we are in a
product. For, we again consider the intersections with the
bigon. Inessential arcs can be handled as before. An arc that
separates one vertex from the other three is like a vertical arc in
the case of the bigons. An arc that separates the vertices into two
pairs separates a disc in the bigon which gives a compression, and so
we are done in this case.

To find the compression, consider an arc $p$ along the quadrilateral
joining a pair of edges in $S$ (see figure~\ref{F:compress}). This is
homotopic to an essentially unique arc in $S$ before puncturing (as we
are in a product). 

We claim that the arc is unknotted, i.e., it is isotopic to this arc
and unlinked from the punctures on which the vertices lie. This
follows since the quadrilateral corresponds to a cancellation, and
hence we can find an embedded, once-punctured disc that the arc $p$
together with an arc in $S$ bounds. Namely, we can join the end points
of this arc by an arc in $S$ running through the bridges that have
been cut and the common island in which these end. The original arc
together with this bound a tube, one half of which is in the chopped
off portion of the plate and the other in the beam corresponding to
the cancellation. The other boundary of this tube is a circle in $S$,
whose boundary is an arc along the bridge that forms a loop and the
other in the other island that our beam bounds. 

If the total space is $F\times I$ with $F$ a punctured sphere, as the
arc is unknotted and unlinked from the punctures containing its
vertices, the arc bounds a compressing disc. In the case of a general
product, there exists a family of annuli, which we can assume are
disjoint from the arc and the circle in $S$ that are part of the
boundary of the quadrilateral, so that the complement of these is
$F\times I$, with $F$ a punctured sphere. By considering outermost
components of intersection as before, we reduce to the previous case.

Thus, we have reduced to the case where we cut along bigons, which is
the previous lemma.
\end{proof}

 We conclude that some handle-structure corresponding to the
presentation obtained by applying the Andrews-Curtis moves is
sphere-like, i.e., the presentation is sphere-like.
\end{proof}

To summarise,

\begin{theorem}\label{T:alg} There is an algorithm that recognises a class of
presentations of the trivial group that includes all those obtained
from the standard one by Andrews-Curtis moves.
\end{theorem}

An immediate corollary is

\begin{theorem}\label{T:bal} At least one of the following holds
\begin{itemize}
\item There is an algorithm to recognise balanced presentations of the
trivial group, or
\item The (balanced) Andrews-Curtis conjecture is false.
\end{itemize}
\end{theorem}

One may use the algorithm given here to attempt to find particular
presentations that violate the (balanced) Andrews-Curtis
conjecture. To do so, we take a candidate presentation of the trivial
group and apply the algorithm to it. If it turns out not to be
sphere-like we have a counterexample.
	
In particular, this method could show that the Andrews-Curtis
conjecture (in the balanced case) is false. Furthermore, if this
method works, for instance, for presentations corresponding to the
Akbulut-Kirby $4$-sphere, we would be able to conclude that property-R
for links is false.

\end{document}